\documentclass[12pt,reqno]{article}

\usepackage{fullpage}
\usepackage{amssymb}
\usepackage{amsmath}
\usepackage{amscd}
\usepackage{mathtools}
\usepackage{float}
\usepackage{amsthm} 
\usepackage{amsfonts}

\begin{document}

\def\Dat{15a}
\def\Version{{\tt 
Klaus Weise 04.10.2021 arx-dynsyst-\Dat.tex 
\Dat.04}}
\def\Lrm{\scriptsize}
\def\Mk{$\spadesuit$ }
\def\En{\enskip}
\def\Eq{\En;\quad}
\def\Ni{\par\noindent}
\def\Mn{\medskip\noindent}
\def\Pn{\bigskip\noindent}
\def\eps{\varepsilon}
\def\ovl{\overline}
\def\unl{\underline}
\def\wit{\widetilde}
\def\Page{\vfill\eject}
\def\fC{f_\text{C}}
\def\fP{f_\text{P}}
\def\fS{f_\text{S}}
\def\Knot{$\bigcirc\!\!\!\to$}
\def\Set{\mathbb}
\def\Ra{\Rightarrow}

\def\be{\begin{equation}}
\def\ee{\end{equation}}

\hskip -.6cm
\Version
\vskip 1cm

\begin{center}
{\LARGE\bf Convergence Criteria for}
\vskip .2cm
{\LARGE\bf Dynamic Integer Systems}
\vskip .6cm
{\large\bf Klaus Weise}
\vskip .4cm
{\large Parkstr. 11, D-38179 Schw\"ulper,
Germany}
\vskip .4cm
{\large\tt kl\underline{~}weise\/@\/t-online.de}
\end{center}
\vskip 1cm

\begin{abstract}
Criteria are presented for testing 
whether every trajectory (integer sequence) 
of a dynamic system converges to the same 
fixed point. The criteria are sufficient and 
applicable to a large variety of systems on  
positive integer values. The variety includes
the (suitably reduced) Collatz 3n+1 dynamic 
system. A system converges according to a 
well-known old criterion if every trajectory value
has a smaller successor. This criterion is
not general enough for application to more 
complex systems. It is therefore generalized
for this purpose by suitable system reduction, 
e.g., by
successively removing branches (values with 
all their own predecessors) without influencing 
the convergence property of the system in 
question until the old criterion becomes applicable.
A system also converges if it has the same 
structure as another system which is already 
known to converge. Moreover, it converges if it
has a self-similar structure where all values 
have either the same or an infinite given 
number of direct predecessors. Every system
is mainly represented and investigated by a 
graph of connected knots identified by the 
trajectory values. The system structure is
then obtained by removing all these values,
but never the knots themselves. The structure 
alone determines the system convergence.
\end{abstract}

\section{Introduction}
\label{Eins}
Natural systems (imagine a fountain or 
a waterfall) develop dynamically and 
irreversibly from a present state to 
following states. A realistic or fictitious 
mathematical model of such a system is 
called a {\it dynamic(al) system\/}. 
Similarly to nature, the model can show 
convergent, divergent, periodic, stable, 
random or chaotic behavior. A dynamic 
system on real values $n$ can be generated 
by iterative application of a unique function 
$f(n)$ to form a sequence of 
{\it successors\/} or {\it values\/} 
$n_{i+1}=f(n_i)$ of a {\it starting value\/} 
$n=n_1$. Values $n_j$ with $j<i$ are {\it 
predecessors\/} of $n_i$. An application 
of the {\it generating function\/} $f(n)$ is 
called a {\it step\/} and the sequence of 
successors is called a {\it trajectory\/} 
$T(n)=(n_1=n,n_2,\ldots,n_i,\ldots)$. 
A trajectory is also called an {\it integer 
sequence\/} if all values $n_i$ are integers.
The aim of the present study is to construct 
sufficient criteria for stating that every
trajectory of a dynamic system on a set 
$\Set V$ of {\it admitted\/} integers $n>0$ 
{\it converges\/} to the same reoccurring 
(or stopping) value $\wit n=f(\wit n)$
called the {\it fixed point\/}.

An outstanding example is the (original) 
{\it Collatz dynamic system} which is generated 
on all integers $n>0$ by the function
\be
\fC(n)=\begin{cases}
~3n+1\En, & \text{if $n$ is odd~;}\\
~n/2\En,  & \text{if $n$ is even~.}
\end{cases}\label{OneA}
\ee
The famous {\it Collatz conjecture\/} 
(which was posed in 1937 and is also 
called the 3$n$+1 conjecture) refers to 
the observed, but seemingly yet unproven 
property of the Collatz system that every 
trajectory leads to a periodically 
reoccurring value 1. The function $\fC(n)$ 
according to Equation (\ref{OneA}) can be
simplified by excluding the needless and 
hindering values $n$ divisible by 2 or 3
(see Section~\ref{Sieben}).
The result of this  reduction is
\be
n'=\fC(n)=(3n+1)/2^\nu \En,
\quad\text{if $n$ and $n'$ 
are indivisible by 2 and 3}
\label{OneB}
\ee
where the integer $\nu>0$ is chosen 
for each $n$ so that $n'=\fC(n)$ 
becomes odd. Then, all successors 
of the value~1 are~1 since $\fC(1)=1$. This 
fact means $T(1)=(1,1,\ldots)$ and is also 
taken as a trajectory stop. The reoccurring 
(or stopping) value $\wit n=1$ is the fixed 
point where all trajectories are conjectured 
to stop.

In Section~\ref{Zwei}, a general {\it dynamic 
(integer) system\/} is introduced and represented 
by a {\it graph\/} of connected {\it knots\/}.
These knots are either arbitrarily, but uniquely 
identified or not at all. If not, then the
graph is called the system {\it structure\/}.
The structure alone determines the system
convergence. There is a large variety of possibly 
convergent systems (including the Collatz 
system) with the same, {\it isomorphic\/} 
structure. Section~\ref{DreiA} explains why 
knot identifications should be removed to
form the system structure that can finally 
lead to a conjectured convergence.
Section~\ref{Drei} deals with system
isomorphism and transformation needed for
using the {\it convergence criteria\/} defined
in Section~\ref{Vier}. These criteria
are then applied in Sections~\ref{Fuenf} 
to~\ref{Sieben} to some systems which are 
already known or not known to converge. 
In particular, the Collatz system turns out 
to converge (see Section~\ref{Sieben}).

There is already an ``old" convergence 
criterion. Accordingly, a system converges 
to the fixed point $\wit n=1$ if every value 
$n>1$ has a smaller successor. Then, this 
successor can repeatedly represent and 
replace $n$ until $n=\wit n=1$ is obtained. 
Another approach could show that 
{\it cycles\/} or {\it divergent\/} 
trajectories (see Section~\ref{Zwei}) do 
not exist, but it was tried in vain to the Collatz 
conjecture. Therefore, essentially other, 
much more general criteria are presented. 
Then, a system converges if it has, 
e.g., a self-similar structure
(see Sections~\ref{Zwei} and~\ref{Vier}).

In this study and if not otherwise stated,
the term ``number" mainly means a result 
of counting or marking objects, whereas
a number itself is called by a more
informative term like ``integer",
``prime" or ``real value". 
This use avoids phrases
like ``... a number of numbers~...".
A lower case italic Latin or Greek letter
always denotes an element of the
set~$\Set N^+$ of all positive integers.
Every used set or infinity means a 
countable one. Phrases similar to 
``finite or (countably) infinite" or
``set or subset" are abbreviated by
``(in)finite" or ``(sub)set", respectively,
a reference to the Collatz case by ``CC".

For more details about approaches to 
the Collatz problem (CC), see, e.g., 
References \cite{Wir} to \cite{Opf} and the 
literature cited there. Wirsching's book 
\cite{Wir} is a comprehensive overview. 
The articles in Wikipedia \cite{Wik} and by 
P\"oppe \cite{Ppp} (in German) are first 
introductions to the matter. Feinstein 
\cite{Fei} and Opfer \cite{Opf} already 
presented interesting, but rather complex 
and seemingly unconfirmed proofs stating
that the Collatz conjecture 
``is not provable" and ``holds true", 
respectively (see Section~\ref{Acht}).

\section{Dynamic system, its graph and structure}
\label{Zwei}
A dynamic (integer) system is usually generated,
represented and investigated by integers $n'=f(n)$.
But this approach can be complemented impressively 
by knots of a graph introduced as follows.

As well known, a (countable) {\it set\/} contains
a counted number 0, 1, 2, ..., or $\infty$ of 
{\it elements\/}. An element is any object or 
may itself be a (sub)set.

Similarly, a (countable) {\it graph\/} is a set of 
a counted number 0, 1, 2, ..., or $\infty$ of 
{\it knots\/}. A~knot is any object or 
may itself be a (sub)graph. The knots are
(dis)connected and uniquely identified 
or not. A  graph of unidentified knots 
is called a {\it structure\/}.

The structure is a
very important essential property of a
dynamic system (say, the system skeleton).
It is in general defined and represented by
the system graph. Such a graph consists of an
(in)finite number of identical knots which are 
arbitrarily (in)directly connected with one 
another or not, and (only for the structure) 
without 
any knot identification by numbers, names or 
other marks. It is stressed that the graph and 
its structure have to be understood purely 
topologically without any dimensions like
space, time, measure, dynamics.

In the following, only dynamic systems 
are considered where each of them is
generated by iterative application of a 
unique function $n'=f(n)$ on positive
integers $n,n'\in\Set V$.
This generating function may not only
be given by arithmetic formulas, but also by 
an algorithm or a list of connections $n\to n'$. 
The system structure is represented 
by a {\it directed graph\/} of connected, 
identical knots \Knot, each with a single 
connection pointer to only one following 
knot, the {\it direct successor\/}. As 
an exception, a knot $\bigcirc$ without 
a pointer may be used for representing 
a fixed point. 
Although knot numbering could diminish 
the transparency of the system structure, 
the values $n$ of the dynamic system 
are later on used as knot numbers
(see Section~\ref{Drei}).

Trajectories of a dynamic system
are {\it connected\/} if they have 
common values. If all trajectories 
are connected with one another, then 
the entire system is connected. A value 
$\wit n=f(\wit n)$ is a fixed point. 
If all trajectories converge to 
the same fixed point, then the system 
is {\it convergent\/}. A periodically 
reoccurring sub-sequence of values of 
a trajectory is called a {\it cycle\/}. 
A fixed point could be taken as a cycle 
with period~1, but this is avoided 
by the present study. Every trajectory 
ends at its {\it root\/}, that is, it 
converges either to a fixed point or to a 
cycle or it {\it diverges\/} to infinity. 
Every root characterizes a {\it subsystem\/}
disconnected from others. Every value can 
act as a {\it delegate\/} (even as the 
fixed point) of all its predecessors. Only
these predecessors together with their 
delegate form a {\it branch\/}. The
delegate is the only root of a branch.
This fact will become important in
Section~\ref{Drei}.

As already described in Section~\ref{Eins}, 
the aim of the present study 
is to introduce criteria for stating that 
an entire dynamic system, generated by 
a function $f(n)$, is connected and 
converges to a single fixed point 
$\wit n=f(\wit n)$. Then, more roots 
do not exist. Before a convergence 
criterion can be applied, the system 
graph has to be constructed by using
the function $f(n)$ with $n$ acting as
the knot identifying 
numbers. A following removal
of all these numbers cleans the graph, 
but does not change the system structure
at all. The structure alone determines 
the roots of the system and should therefore 
be investigated, whereas $f(n)$ itself with 
all values and trajectories will play a minor 
role only.

It is stressed again that a structure 
only deals with knots and their 
connections, but yet not with the knot 
identification, e.g., by values $n$
of the dynamic system in question.
Systems with the same structure are 
{\it isomorphic\/}. A set of them is 
called a {\it family\/}. Two structures 
are isomorphic if there is a complete 
one-to-one (bijective) 
correspondence between them, 
concerning all knots (unidentified),
connections, roots, trajectories.

There are related structures (typed by 
$\eta$) which can be used as valuable 
tools in convergence investigation.
Every chosen one of these structures 
defines an isomorphic and convergent 
infinite system family. And every knot X of 
the chosen structure points to its own direct 
successor and has a given, for all knots
the same fixed counted number $\eta>0$ or 
$\eta=\infty$ of other knots which are the 
{\it direct predecessors\/} of knot X. Only 
a single knot does not have a successor, 
the fixed point of the chosen structure.
All knots are (in)directly
connected with one another at least
at the fixed point. Therefore, other
roots (fixed points, cycles and 
divergent trajectories) do not exist. 
 
{\it Self-similarity\/} is another 
important essential property of every 
graph (or system) of a family with 
the just chosen particular structure. 
One can easily see by inspecting the
graph that every knot acts as the fixed 
point of a branch and every branch 
has exactly the same structure as the 
chosen structure of the entire graph. 
This property is referred to
in the following by the adjective
``self-similar (of type~$\eta$)".

More generally, only infinite, connected 
dynamic systems with a single fixed point 
and no other roots can be self-similar  
provided that the fixed point and every 
knot have the same number $\eta>0$ 
or $\eta=\infty$ of direct predecessors.
Otherwise, the structure of the entire
system and that one of every branch cannot
be isomorphic (see Sections~\ref{Fuenf} 
to~\ref{Sieben} for examples).

\section{Interchanging or removing knot identifications}
\label{DreiA}
The following story shall explain that knots of
a system graph do not necessarily need to be 
identified and available to determine system 
roots, e.g., convergence to a fixed point.
Any marks used to identify knots can 
arbitrarily be interchanged or removed without 
modifying the system structure with all its
knots, connections, roots.

Consider a town with many sights and 
other important objects such as places, 
large and small streets, crooked lanes, 
churches, monuments and more. Every 
day, many tourists arrive at the railway 
station and want to return in the 
evening with the last train. A problem 
is that many tourists get lost in the 
town and cannot be back in time at the 
hidden station. And there are not enough 
hotels. What should be done? Two 
alternative solutions are discussed
by the town council:\Ni
(a) Every 
tourist gets a detailed town map when he 
starts at the station for sight-seeing.\Ni
(b) At every important object of the town 
there is a signpost ``To the station $\to$". 

Walking to the station can be taken as a
dynamic system with the station as a fixed 
point and the important objects of the town
as the values or knots. The Solutions (a) and 
(b) are equivalent alternatives.
With Solution (a), every object must be uniquely
identified (e.g., by a name) itself and on
the map as well. Solution (b) represents the 
structure of the system, a directed graph of
the object connections. Object identification  
is not needed, but it must be required that
a pointed way (trajectory) does not end at a 
cycle or anywhere off the station. This can 
best be achieved by setting up the signposts 
(knots) from the station back to the objects. 
There is no need to identify the signposts,
but if they are identified (by names or numbers 
for maintenance or other purposes), then
they can arbitrarily be interchanged or not 
without any modifying the structure. This 
sentence is most essential. It expresses a new
important finding stated here in general for short: 
Structure only determines convergence and 
roots at all !

\section{System isomorphism and transformation}
\label{Drei} 
As just stated above,
only the structure can determine
convergence of a dynamic system~A to a single
fixed point. System~A is proven to converge if
it has the same structure as another system~B 
which is already known to converge 
(see Section~\ref{Vier} Criterion~1). 
This {\it isomorphism\/} of A and~B is often 
missing, but can possibly achieved by system 
transformation which is {\it equivalent\/} with
respect to the main intention, e.g., to prove a 
conjecture or to reduce or to expand a system.

Pay first attention to several important general 
building blocks for equivalently transforming 
(reducing or, reversely, expanding as well)
a dynamic system which shall be tested for
convergence and is represented by a
directed graph of knots \Knot ~and $\bigcirc$.
Let these knots be uniquely identified by
admitted numbers $n\in\Set V$ and each knot 
be connected and pointing at most to one 
single direct successor knot $n'\in\Set V$ 
obtained by $n'=f(n)$. The main question is
whether the system converges to a single
fixed point $\bigcirc$ $\wit n=f(\wit n)$.
This equation must have just one single 
solution $\wit n\in\Set V$. Otherwise, the 
system does not converge. The system 
structure is formed by removing all knot 
identifying numbers. But all knots with 
all their directed connections remain.

Every knot can not only be
arbitrarily and uniquely identified 
by a number $n$ or otherwise or not at all.
It can itself, if suitable, also be a connected 
(sub)graph. Reversely, any connected 
(sub)graph can also act as a single knot. 
This slightly generalized knot definition 
easily allows equivalent
system transformations, preferably 
reductions with respect to convergence
or, as well, expansions to make systems 
isomorphic.

Notice that such a transformation changes
not only the graph and the structure,
but also the generating function $f(n)$ and
the set $\Set V$ of admitted knot numbers.

Examples of transformation building blocks
for direct or reversal application:\Ni 
{\bf Block (1):~} A chain \Knot\Knot...\Knot
~with all its {\it inputs\/} (connections from 
the direct predecessors) can be replaced by 
a single knot \Knot~with all inputs 
from outside the chain.\Ni
{\bf Block (2):~} A cycle can be replaced 
by a fixed point $\bigcirc$ with all inputs
from outside the cycle.\Ni
{\bf Block (3):~} A knot \Knot~with no 
input at all can be removed.\Ni
{\bf Block (4):~} A branch $B(n)$ (knot $n$
\Knot~together with all its predecessors, see
also Section~\ref{Zwei})
can be taken as a subsystem converging 
to knot $n$ and can completely be removed
since the branch has only the single root
knot $n$ (delegate).
But attention! Knot $n$ itself
must not be a fixed point, a cycle member, 
or infinite. Otherwise, the due (possibly
questionable) subsystem would disappear 
and could thus cause a wrong proof of 
convergence. If $n'>n$, it can sometimes
be sufficient to interchange the knot 
numbers $n$ and $n'$ only instead of
removing the whole branch $B(n)$.\Ni
{\bf Block (5):~} Knot numbers, names, 
or other knot identifiers can arbitrarily be 
interchanged (by permutation), replaced 
or removed without any modifying the 
structure\Ni
{\bf Block (6):~} If two trajectories are
connected at some knot, then an arbitrarily
other knot of both trajectories can instead
be chosen for the connection.

A short example of using the building blocks 
is the equivalent transformation of CC from 
Equation~(\ref{OneA}) to Equation 
(\ref{OneB}). Blocks (1) and (3) are used
to exclude values $n$ divisible by 2 and~3,
respectively. For a more elaborate 
transformation, see CC in Section~\ref{Sieben}.

Let now again an (in)finite set $\Set V$ 
be given, the elements of which are all 
the admitted values $n$ 
of a dynamic system generated by the 
unique function $f(n)$. The values $n$
shall serve for uniquely identifying the 
knots of the system graph. If a person~A 
needs all $n\in\Set V$ for the identification, 
then a person~B also needs all $n\in\Set V$
provided that these elements $n$ differ 
from one another and each $n$ is used 
once and only once. Naturally, the 
distributions of the values $n$ to 
the knots by person~A and by person~B
differ in general by a permutation.

Consider again a system as just described with
an (in)finite number of knots. Let every knot
be uniquely identified by a value $n\in\Set V$
and have an own given (in)finite number
$\eta(n)\ge0$ or $\eta(n)=\infty$ of direct 
predecessors forming a set $\Set U_n$. 
All these $\Set U_n$ are disjoint. This fact 
allows altogether a unique serial numbering 
of all knots since it is well-known that the 
union of (in)finitely many (in)finite, disjoint 
sets $\Set U_n$ becomes a single (in)finite 
set $\Set U=\bigcup_n^{\eta(n)}\Set U_n$. 
It is stated again that all the 
admitted $n\in\Set V$ have to be
different and used once and only once 
for uniquely identifying the knots.

Although knots are usually numbered
by iterative application of the unique 
generating function $n'=f(n)$, there is also 
another, possibly more practical way to 
generate numbers for the knots of a 
conjectured convergent structure. This way 
is an iterative application of the reverse
generating function $n=f^\bullet(n')$ with
start at the already known fixed point
$n'=\wit n$. The reverse function is not
unique. Every $n'$ can have an own (in)finite
number of results $n$, the direct predecessors 
of $n'$. These different results $n$ 
should also be identified themselves by a 
suitable new parameter $\mu$. The described 
way identifies all knots of the convergent 
subsystem with $\wit n$ as its root. If all 
admitted $n\in\Set V$ are needed, then no 
$n$ at all remains for other, disconnected 
subsystems and, thus, the entire system 
is proven to converge. 

\section{Convergence criteria}
\label{Vier}
After laying the foundations in
Sections~\ref{Eins} to \ref{Drei}, 
six convergence criteria 
can now easily be constructed and 
understood. In every system, all the
admitted $n\in\Set V$ have to be used 
once and only once for a unique knot
identification. The criteria are 
mainly based on the idea of strictly 
separating structure from knot identification.
For applications of the criteria, see  
Sections~\ref{Fuenf} to~\ref{Sieben}. 
The criteria are to be used to test 
the convergence of a dynamic system to 
the fixed point $\wit n$. The criteria are
sufficient. This means that in case of
a negative test result, one cannot decide 
whether or not the system converges.
Before a criterion is applied, the system 
should be suitably reduced, e.g., by 
using the transformation building blocks 
of Section~\ref{Drei}.
The following criteria together with
system transformations can turn out to be
related to one another, e.g., Criteria 3 and 4.
It is also stressed again that the system
to be tested for convergence must have 
one single fixed point only.

\Mn{\bf Criterion 1:~}
If two dynamic systems A and B 
are isomorphic (have the same structure)
and system B is known to converge, 
then system A also converges.\Ni
--- Comment: Structure alone 
determines roots 
(see Sections~\ref{Zwei} and \ref{DreiA}).

This convergence Criterion 1 is quite simple and
clear, most basic, very general, and does not need
any knot identification. But a convergent system~B
must be available for the structure comparison. 
The criterion could even be more generalized for
other structural system properties, e.g.,
cycles or divergences.

The systems~A and B 
are isomorphic if there is a complete 
one-to-one correspondence between 
them, concerning all knots (unidentified),
connections, roots, trajectories.
If A and B are not isomorphic, then
they can possibly equivalently transformed to 
become isomorphic, e.g., if the structure 
of A (or~B) contains that of B (or~A) with 
corresponding fixed points. In particular, this 
applies to the case of a self-similar A where the
entire A and every branch are isomorphic.
The branch converges and can act as~B. The 
best way to test A and B for being isomorphic 
seems trying to generate the graph of A step by 
step with the reverse function $n=f^\bullet(n')$
on the known graph of B from the 
known fixed points $\wit n$ of A and B up 
to all predecessors of every knot $n'$ of A.
Then, system A converges if all $n\in\Set V$ 
are needed. (See also Section~\ref{Drei}
and Criterion 2.)

\Mn{\bf Criterion 2:~}
A system converges if all its knots are (in)directly 
connected with one another and need all
$n\in\Set V$ for a unique identification.\Ni
--- Comment: No disconnected knots 
remain for trajectories to roots besides the
single fixed point (For application, see the
last paragraph of Section~\ref{Drei}
and also Criterion 1).

\Mn{\bf Criterion 3:~} 
An infinite system converges 
if its structure is self-similar. \Ni
--- Comment: Systems which are finite
or have roots besides the single fixed
point are never self-similar 
(see end of Section~\ref{Zwei}).

\Mn{\bf Criterion 4:~} 
An infinite system converges if all knots 
have the same number $\eta>0$ or
 $\eta=\infty$ of direct predecessors.\Ni
--- Comment: The system is self-similar 
(see Criterion~3).
Moreover, systems with the same 
$\eta$ of all knots are always isomorphic.

\Mn{\bf Criterion 5:~} 
A system converges if it completely consists 
of branches with all their delegates (in)directly 
connected with one another.\Ni
--- Comment: The branches can successively be
removed without influencing the convergence
to the fixed point
(see Section~\ref{Drei} Blocks (1),(3),(4)).

\Mn{\bf Criterion 6:~} 
A system converges to the single fixed point 
$\wit n=\min n\in\Set V$ if every knot 
$n>\wit n$ has a successor knot $m<n$.\Ni
--- Comment: This is the old criterion 
(see Section~\ref{Eins}). 

The old criterion can only be applied
if every case $n'=f(n)>n$ can be removed by a 
suitable equivalent system transformation, e.g., by
removing the branch $B(n)$ or by interchanging
the knot numbers $n$ and $n'$ (only these 
numbers, never the knots themselves with 
their directed connection!).
Cases $n'>n$ of cycles or divergent trajectories
can never be removed at all. (See also
Section~\ref{Drei} Blocks (4) and (5).) But
such cases need no attention. This is an
advantage of the old criterion. A shortcoming
seems to be that structure is not taken into
account.

\section{A simple dynamic system}
\label{Fuenf}
The convergence criteria are first 
applied to the simple (but not trivial) 
dynamic system with the following 
generating function $\fS(n)$ and its 
reverse $\fS^\bullet(n')$:
\be
n'=\fS(n)=\begin{cases}
~n\En,         & \text{if $n=1$~;}\\
~(n-1)/2\En, & \text{if $n>1$ is odd~;}\\
~n/2\En,      & \text{if $n$ is even~;}
\end{cases}
\label{SimA}
\ee
\be
n=\fS^\bullet(n')=2n'\En \text{and~} =2n'+1
\label{SimB}
\ee
with all values $n,n'>0$ admitted and 
forming the set $\Set V=\Set N^+$. Since 
always $n'<n$, except $n'=n=\wit n=1$ 
for the fixed point $\wit n$, the system 
is connected and converges to $\wit n$ 
according to the old Criterion 6. 
Every $n'$ has  
$\eta=2$ direct predecessors. This fact 
makes the system structure self-similar. 
Accordingly, the system also turns out 
to converge by using Criteria 3 or 4.

For a test of convergence to $\wit n$,
the knots of a system graph can best 
also be generated and numbered by
$n=f^\bullet(n')$ reversely from $\wit n$.
If all $n\in\Set V$ are needed once and 
only once, then the system converges
since no values $n$ remain for 
disconnected subsystems of other roots.
Indeed, $\fS^\bullet(n')$ generates $\eta=2$
knots $n$ for every $n'$ already given. 
Every $n$ is needed 
once and only once and obtains an own 
trajectory $T(n)$ different from others. 
See reversely, e.g.,
$1\Ra2,3\Ra4,5,6,7\Ra
...\Ra2^\nu,...,2^{\nu+1}\!-\!1
\Ra ...$ ($\Ra$ denotes 
a bundle of reverse connections
from ... to ... by $f^\bullet$).
No $n$ is lost and knots with the same 
number do not exist. The system thus converges.

\section{Families of convergent dynamic systems}
\label{Sechs}
The application of the most general convergence
Criterion 1 to a system A requires an available,
possibly isomorphic, convergent system B.
Such a system B can easily be constructed
because every integer $n>1$ can uniquely be 
represented by multiplied primes or, similarly, 
by a sum of powers $2^\nu$ with exponents 
$\nu\in\Set N$, or by proceeding one of
many other ways. 

In particular, the generating function $\fP(n)$ 
of a large family of isomorphic systems on 
multiplied primes (MP) and their reverse 
function $\fP^\bullet(n')$ read
\be
n'=\fP(n)=n/p(n) \Eq
n=\fP^\bullet(n')=n'p \Eq
(\En p, p(n)>3\En)
\label{Prim}
\ee
where $p$ and $p(n)$ are arbitrarily 
chosen prime factors (of $n$ for $p(n)$).
The particular choice $p, p(n)>3$ of 
Equation~(\ref{Prim}) is suitable only for
use in Section~\ref{Sieben}. The systems 
of the family differ only by the sequences
of chosen primes used for generating 
trajectories, that is, by knot numbering. 
Accordingly, the systems are isomorphic
(all have the same structure). Since always 
$n'<n$ and by application of Criterion 6, 
all the systems converge to the only 
fixed point $\wit n=1$. A more detailed
example follows.

The knot numbers $n$ of the family graphs
are obtained by $n=\fP^\bullet(n')$ with an
infinity $\eta=\infty$ of direct predecessors 
$n$ of every given $n'$, beginning with 
$n'=\wit n=1$. For example with all $p>3$, 
the direct predecessors of 
the fixed point $n'=\wit n=1$ are all single
primes $n=p_1=5,7,11,...$ in arbitrary 
order; those ones of $n'=5$ are all 
the products $n=p_1p_2=25,35,...$ 
of two primes, also in arbitrary
order, and so on. Every 
direct predecessor $n$ of $n'$ has one prime 
factor more than $n'$ and is indivisible 
by 2 and 3 so that $n=6k\pm1>1$ for every
integer $k>0$ (see also CC in Section~\ref{Sieben}).
All these $n$ and $\wit n=1$ 
form the infinite set $\Set V$ of admitted
values. Every $n$ has to be used as a knot 
number, but only once for a unique knot
numbering. Indeed, all direct predecessors 
$n$ of some given $n'$ differ from one 
another due to different prime factors $p$. 
They also differ from all other values $n$ 
because of a different number of prime factors. 
And every knot $n$ can be arrived by a 
reverse trajectory from $\wit n=1$ successively 
multiplied with all prime factors of $n$. 

Quite similarly, the generating function 
$f_2(n)$ of another large system family 
on sums of powers $2^\nu$ and the
reverse $f_2^\bullet(n')$ read
\be
n'=f_2(n)=n-2^\nu \Eq
n=f_2^\bullet(n')=n'+2^\mu
\label{Sum}
\ee
where $\nu$ is a suitably chosen exponent 
and $n>\wit n=0$. All systems converge
according to Criterion 6 since always $n'<n$.
Again, every knot $n$ can be arrived by a
reverse trajectory from $\wit n=0$ successively 
added with all powers $2^\mu$ of $n$.

All systems of both and related other 
families are self-similar since $\eta=\infty$ 
of every knot. They also converge 
according to Criteria 3 or~4.
Moreover, the systems are isomorphic, i.e.,
all have the same structure independent
from any unique knot numbering. 

\section{Collatz dynamic system and conjecture}
\label{Sieben}
The convergence criteria of Section~\ref{Vier}
together with the graphic properties, tools 
and transformations of Sections~\ref{Zwei} 
and~\ref{Drei} are now applied to easily 
prove the convergence of a suitably reduced 
Collatz system for confirming the Collatz 
conjecture (CC, see Section~\ref{Eins}).

The original CC system graph is first generated
by $n'=\fC(n)$ of Equation~(\ref{OneA}) on all
admitted knot numbers $n,n'>0$. The graph 
has a cycle $T(1)=(1,4,2,1,\ldots)$ with $n=1$ 
conjectured to be also a value of every 
trajectory $T(n)$. The graph is then reduced 
suitably and equivalently with respect to the
CC conjecture by using the building blocks
of Section~\ref{Drei} for transformations 
until a criterion can work. 

{\bf 1st reduction:~}
According to Block~(1), every chain of knots, 
generated by some even $m$ divided $\nu$ times
repeatedly by 2, is replaced by the single knot 
with odd $m/2^\nu$. This transformation
removes all knots with even $n$ and leads to
$\fC(n)$ in Equation~(\ref{OneB}) with
$\nu>0$. The cycle disappears, but its member
knot $n=1$ remains as a conjectured fixed 
point $\wit n=1$ of every trajectory. This
fixed point $\wit n=1$ together with 
$\nu=2$ is the only solution of 
$\wit n=\fC(\wit n)=(3\wit n+1)/2^\nu$
transformed to $(2^\nu-3)\wit n=1$.

{\bf 2nd reduction:~}
Every knot with $n$ divisible by 3 has no 
predecessor $m$ and can thus also be 
removed according to Block~(3). No  
$n=(3m+1)/2^\nu$ is divisible by 3.

{\bf 3rd reduction:~}
Block~(4) allows to also remove branches 
$B(n)$, e.g., if $n'>n$ obtained from 
$n'=(3n+1)/2^\nu>n$ for $\nu=1$ only
(see Equation~(\ref{OneB})). Knot $n'$ and 
all its direct predecessors $n$ for $\nu>1$
remain with $n'<n$
since they do not belong to $B(n)$.

Notice that the reductions remove some 
at first admitted $n$ and change the
original cycle to a fixed point, but do not
remove any possibly existing roots.
Similarly, this applies also to the 
branches. The generating 
function $\fC(n)$ of the three times 
reduced CC system graph and its 
reverse $\fC^\bullet(n')$ now read
\be
n'=\fC(n)=(3n+1)/2^\nu \Eq
n=\fC^\bullet(n')= (2^\mu n'-1)/3 
\label{Coll}
\ee
with admitted $n$ and $n'$ indivisible 
by 2 and 3 (or $n,n'=6k\pm1>0$, see below) 
and with the single fixed point $\wit n=1$.
The parameter $\nu>1$ (in contrast to
$\nu>0$ in Equation~(\ref{OneB})!)
serves for always $n'<n$ and makes the old
Criterion 6 thus applicable. Reversely, the 
parameter $\mu>1$ uniquely identifies all 
the infinite number $\eta=\infty$ of 
admitted results (direct predecessors) 
$n=\fC^\bullet(n')$ for every admitted
$n'$ given. Although $\eta=\infty$, by far 
not all values of $\nu$ and $\mu$ have
be used in Equation~(\ref{Coll}) to
connect admitted $n$ and $n'$. But this
fact does not matter.

The CC system is now already proven 
to converge to the fixed point $\wit n=1$
according to the old Criterion 6 since
always $n'<n$. It also converges according
to Criteria 3 and 4 since it is self-similar
because of $\eta=\infty$. See also the
self-similar particular structure chosen in 
Section~\ref{Zwei} and the MP and other 
structures discussed in Section~\ref{Sechs}.
All these structures and the Collatz one are
isomorphic to one another. They are taken
as B and A, respectively, for Criterion 1.
Each of the results here obtained already
confirms the Collatz conjecture.

Nevertheless, someone may remain unsatisfied
possibly by lack of numeric calculations 
similar to those vain ones in the past. More
elaborate arguments shall convince him. 
Read paragraph ``The knot numbers ..." 
of Section~\ref{Sechs} and choose any MP 
system according to Equation~(\ref{Prim}) 
with all the prime factors $p,p(n)>3$.
These factors are needed to generate the 
entire convergent MP graph reversely by
$n=\fP^\bullet(n')$ from the fixed point 
$\wit n=1$ to all knots uniquely numbered 
by using the set $\Set V$ with all the 
admitted $n,n',\wit n$ indivisible by 2~and~3. 
Every $n'$ turns out to have $\eta=\infty$
direct predecessor knots. The entire CC 
graph is now reversely generated by
$n=\fC^\bullet(n')$ strictly in parallel to
MP with the same set $\Set V$. The same 
result $\eta=\infty$ of MP and CC then 
leads with Criteria 3 and 4 to the same
self-similar and convergent structures.
Criterion 1 confirms the convergence 
of CC as system A and MP as system B.
These systems are isomorphic, the
knot number distributions differ only 
by a permutation. This fact becomes clear
since every knot number of CC can as 
well be uniquely expressed by a knot 
number of MP and vice versa. This 
one-to-one correspondence between 
$q=6k\pm1>0$ of CC and ``${\it\Pi}>0$ 
indivisible by 2 and 3" of MP is 
(as required by a criticism) shown 
as follows.

Every $q$ is obviously indivisible by 
2 and 3. But can every ${\it\Pi}$ also uniquely 
be represented by $q$? Here, ${\it\Pi}$ 
is a product of prime factors $p>3$. 
Then, ${\it\Pi}=q=6k\pm1>0$ must apply 
for every ${\it\Pi}$, that is, the unknown
$k$ must have a unique solution.
With other words, either only ${\it\Pi}+1$ or 
only ${\it\Pi}-1$ must be divisible by 6. Both
are already even, divisible by 2. Every
triple of succeeding integers
$a-1, a, a+1$ has always one integer 
divisible by 3. But $a\ne{\it\Pi}$ since 
${\it\Pi}$ is already indivisible by 3. Thus,
either only ${\it\Pi}-1$ or only ${\it\Pi}+1$ is
divisible by 3 and 2. The division by 6 
then leads uniquely to $k$. Accordingly,
CC can be transformed to MP without
modifying the structure and vice versa.
CC and MP are isomorphic. 

This fact can also be stated more easily since
every knot $n'$ has an infinity $\eta=\infty$
of direct predecessors $n$. This results for CC
 from
$n=\fC^\bullet(n')=(2^\nu n'-1)/3$ for an 
infinity of values $\nu$ and all $n$ and $n'$
indivisible by 2 and 3; for MP from
$n=\fP^\bullet(n')=n'p$ for an infinity of
primes $p>3$ and all $n$ and $n'$. But these 
values do not matter at all. They only show 
a one-to-one correspondence between
the knots of CC and MP. 
Much more easily, if each of all knots 
of two systems A and B has the same number
$\eta>0$ or $\eta=\infty$ of direct 
predecessors, then A and B are isomorphic.

\section{Concluding remarks}
\label{Acht}
As the main result of the present study,
sufficient criteria are established for testing
whether a dynamic system on positive integers 
is connected and converges to a fixed point. 
The criteria are based on a quite simple idea. 
Let, e.g., the structure of a
dynamic system in question be represented 
by an infinite, directed graph of identical knots 
\Knot, each with a single connection pointer to 
only one following knot. If this structure is 
self-similar or the same as that of another 
system already known to be convergent,
then cycles and divergent
trajectories do not exist and the system 
converges to just a single fixed point.
Every knot identification, e.g.,
by numbers does not influence the
system structure. The criteria thus are
general enough for application to a large
variety of related systems and problems. 
In particular, the Collatz conjecture
is easily confirmed. The criteria also
allow to avoid knot 
numbers at all to easily overcome the 
high barrier between very many logical
steps and an infinity of them.

One may ask why the properties of the huge 
number of previously investigated individual 
trajectories of the Collatz dynamic system (see 
Wikipedia \cite{Wik} or P\"oppe \cite{Ppp}) were 
not taken into account. The reason is that the 
present approach only considers the structure 
and one characteristic problem of the system 
in its entirety, namely, its convergence to the 
fixed point~1. Consider, e.g., the related system 
generated by the function $f(n)=(n+1)/2^\nu$
($n$ and $f$ are odd).
Here, every trajectory is proven to converge 
to the fixed point~1, easily and merely by 
the general and inherent property $n'=f(n)<n$ 
of $f(n)$ itself $(n>1)$, but not by taking 
into account infinities of trajectories.

Many experienced mathematicians on number 
theory, inspired laymen, and, previously 
also this author, tried in vain to prove the 
Collatz conjecture. Why did they not succeed? 
Possibly, they followed a mainstream approach 
in which too much attention was paid to all 
of the trajectories (e.g., by relying 
on computer experiments) rather than paying 
attention to the problem in its entirety, say,
its infinite structure.

Similarly, Feinstein \cite{Fei} may be right in
that an infinity of program lines or computing 
time is needed to test whether or not all the 
individual trajectories converge to~1. However, 
to solve a problem it may be sufficient to 
investigate an essential property, such as 
the system structure or
a permanent increase of entropy (or loss of 
information), which allows for a completely 
different approach that avoids difficulties. 
One example is the irrational number 
$\pi$ that can never be determined exactly by 
numerical calculations, but many essential 
facts about $\pi$ are already known and always 
used. The human brain can think better than 
a computer about abstractions, continuities, 
infinities and irrationals. Therefore, 
there are not necessarily contradictions 
between Feinstein's proof \cite{Fei}, Opfer's 
proof \cite{Opf} and the present approach. 
Problems can have several solutions, 
impractical ones needing an infinity of 
logical steps or others requiring a large 
or short finite number of them. See, for 
instance, $\Sigma_{i=0}^\infty x^i=1/(1-x)$
with $|x|<1$.
Then, one can say that the present effort 
to confirm the Collatz conjecture is short 
compared with other approaches and could 
finally finish the research of more than 
80 years.

\section{Acknowledgments}
\label{Neun}
The present paper describes results of the 
author's private theoretical research. Some applied 
concepts, facts, and methods are well-known 
in mathematics, especially in number theory 
of dynamic systems. They are combined or
used as tools to establish new, simple
convergence criteria for dynamic systems.

The author would like to express his thanks to 
the Physikalisch-Technische Bundes\-anstalt (PTB, 
Federal Institute of Physics and Technology,
Braunschweig, Germany) for years of support, 
to Prof. Dr. L. Baringhaus (retired from 
Leibniz University of Hannover, Germany), Dr. 
M. Matzke (like the author retired from PTB), 
and Dr. M. Reginatto (PTB) for suggestions and 
criticisms concerning the convergence problem 
and also for their patience with the author's 
several vain trials to solve the problem. 
He also thanks Dr. M. Matzke, Prof. Dr. R.
Michel (also retired from Leibniz University) 
and Dr. W. W\"oger (also retired from PTB) 
for decades of successful collaboration in 
theoretical physics of ionizing radiation and 
in general measurement, Prof. Dr. L. 
Baringhaus for valuable advice in stochastics, and 
P. Harris and the author's granddaughter S. Weise
for linguistically inspecting the English text.

\section{A motto for dynamic systems}
\label{Zehn}
The following poem, ``The Roman Fountain" by 
Conrad Ferdinand Meyer (Swiss novelist and 
epic poet, 1825--1898), could be a motto for 
the dynamic systems, each appearing as a 
fountain of integers.

\centerline{\bf Der r\"omische Brunnen}

{\obeylines\everypar{\hfil}\parindent=0pt           
Auf steigt der Strahl und fallend gie\ss t
Er voll der Marmorschale Rund,
Die, sich verschleiernd, \"uberflie\ss t
In einer zweiten Schale Grund;
Die zweite gibt, sie wird zu reich,
Der dritten wallend ihre Flut,
Und jede nimmt und gibt zugleich
Und str\"omt und ruht.}

\Pn 
{\large\bf The Author}

\Pn
Prof. Dr. Klaus Weise\Ni
Parkstr. 11\Ni
D-38179 Schw\"ulper, Germany\Ni
E-mail:~ kl{\tt\underline{~}}weise@t-online.de\Ni
\En\En(Please note the underscore after kl)\Ni 
Phone:~ +49-(0)5303-5108\Ni
Files:~ arx-dynsyst-\Dat.tex and .pdf

\Pn
The author, Dr. rer. nat. Klaus Weise (born 
in 1934), is a theoretical physicist in the 
fields of general measurement and ionizing 
radiation. He worked as Director 
and Professor at the Physikalisch-Technische 
Bundesanstalt (PTB, Federal Institute of 
Physics and Technology) in Braunschweig, 
Germany. He also was lecturer in theoretical 
physics at the Leibniz University of 
Hannover, Germany. He is (main) author of 
numerous research publications and several 
(inter)national standards. He retired 
from his working places but continued 
his activities in his research fields.
As a hobby, he tries to solve unproven
mathematical problems.

\Pn\Version

\end{document}